\documentclass[12pt, letterpaper]{article}

\usepackage{amssymb}
\usepackage{amsfonts, amsmath}

\begin{document}
\begin{center}
\textbf{The Extended Generating Function of the Radical of $\mathbf{n}$ \\ and the abc-Conjecture}\\ \vspace{0.5cm}
\textbf{Constantin M. Petridi}
\\ cpetridi@math.uoa.gr
\\ cpetridi@hotmail.com
\end{center}
\vspace{1cm}
\textbf{Abstract}. We introduce the function defined by the sum of the generating series $\sum_{n\geq{1}}\dfrac{R(n)^{t}}{n^{s}}$, where $R(n)$ is the radical of $n$ and $s$ and $t$ real positive variables. The difference with the ordinary generating series $\sum_{n\geq{1}}\dfrac{R(n)}{n^{s}}$ [1] is that now $R(n)$ appears elevated to the positive power $t$. Since $R(n)$ is multiplicative, logarithmic differentiation with respect to $s$ and $t$ of the series and of its equal Euler product gives an identity involving $s$ and $t$ and two positive functions $\mathcal{S}(s,\,t)$ and $\mathcal{T}(s,\,t)$, expressed as series running over all primes. Appropriately interpreted this identity leads to a proof of Bombieri's $abc$-conjecture/question $a+b<R(abc)^{2}$ for all $n$ satisfying $n<R(n)^{\frac{\mathcal{S}(s,\,t)}{\mathcal{T}(s,\,t)}}$.\vspace{1cm}\\

\hspace{-0.7cm}\textbf{1 The series\,\, $\sum_{n\geq{1}}\dfrac{R(n)^{t}}{n^{s}}$}\vspace{0.5cm}\\
We consider the series $\sum_{n\geq{1}}\frac{R(n)^{t}}{n^{s}}$, where $R(n)$ is the radical of $n$ and $s$ and $t$ real positive variables. Since $R(n)\leq n$ and is only equal to $n$ for the squarefree numbers, it follows that
\begin{equation*}
\sum_{n\geq 1}\dfrac{R(n)^{t}}{n^{s}} < \sum_{n\geq 1}\dfrac{1}{n^{s-t}}.
\end{equation*}
As the series on the right side converges for $t>0$ and $s>1+t$ (region of convergence, denoted for brevity RC) so does the series on the left side. The sum therefore is a well defined function of $s$ and $t$ within RC.\vspace{1cm}\\

\hspace{-0.7cm}\textbf{2 The Euler product of\,\, $\sum_{n\geq 1}\dfrac{R(n)^{t}}{n^s}$}\vspace{0.5cm}\\
\textbf{Theorem 1}.\hspace{0.5cm}
If $s$ and $t$ are within RC, then
\begin{equation*}
\sum_{n\geq 1}\frac{R(n)^{t}}{n^{s}}=\prod_{p}\big(\dfrac{p^s-1+p^t}{p^s-1}\big)
\end{equation*}
where $p$ runs over all primes.\vspace{1cm}\\
\textbf{Proof}. Since $R(n)$ is multiplicative, so is $R(n)^t$ because of $R(1)^t=1$ and $R(nm)^{t}=R(n)^{t}R(m)^{t}$ for coprime integers $n$ and $m$. By applying Euler's generalized identity [1], [2], we have successively
\begin{align}\nonumber
\sum_{n\geq{1}}\frac{R(n)^{t}}{n^{s}}&=\prod_{p}\big(1+R(p)^{t}\frac{1}{p^{s}}+R(p^2)^{t}\frac{1}{p^{2s}}+\cdots\big)\\ \nonumber
&=\prod_{p}\big(1+\frac{p^{t}}{p^{s}}+\frac{p^{t}}{p^{2s}}+\cdots\big)\\ \nonumber
&=\prod_{p}\big(1+\frac{p^{t}}{p^{s}}(1+\frac{1}{p^{s}}+\frac{1}{p^{2s}}+\cdots)\big)\\ \nonumber
&=\prod_{p}\big(1+\frac{p^{t}}{p^{s}}(\frac{1}{1-\frac{1}{p^s}})\big)\\ \nonumber
&=\prod_{p}\big(1+\frac{p^{t}}{p^{s}-1}\big)\\ \label{1}
&=\prod_{p}\big(\frac{p^s-1+p^{t}}{p^{s}-1}\big),
\end{align}
where $p$ runs over all primes. Q.E.D.\vspace{1cm}\\

\hspace{-0.7cm}\textbf{3 Differentiation of\,\, $\ln\sum_{n\geq 1}\dfrac{R(n)^{t}}{n^{s}}=\ln\prod_{p}\big(\dfrac{p^s-1+p^{t}}{p^{s}-1}\big)$}\vspace{0.5cm}\\
Taking the logarithms of both sides of \eqref{1} we get
\[\ln\sum_{n\geq 1}\frac{R(n)^{t}}{n^{s}}=\sum_{p}\ln\big(\frac{p^s-1+p^{t}}{p^{s}-1}\big).\]
Partial differentiation with respect to $s$ of both sides gives
\begin{align}\nonumber
\dfrac{-\sum_{n\geq 1}\dfrac{R(n)^{t}}{n^{s}}\ln n}{\sum_{n\geq 1}\dfrac{R(n)^{t}}{n^{s}}}&=
\sum_{p}\dfrac{p^{s}\ln p(p^s-1)-(p^s-1+p^t)p^s\ln p}{(p^s-1)(p^s-1+p^t)}\\ \nonumber
&=\sum_{p}\dfrac{(p^s-1-p^s+1-p^t) p^s \ln p}{(p^s-1)(p^s-1+p^t)}
\end{align}
and hence
\begin{equation}\label{2}
\dfrac{\sum_{n\geq 1}\dfrac{R(n)^{t}}{n^{s}}\ln n}{\sum_{n\geq 1}\dfrac{R(n)^{t}}{n^{s}}}=
\sum_{p}\dfrac{p^s}{p^s-1}\dfrac{p^t}{p^s-1+p^t}\ln{p}\,\,\,:=\,\mathcal{S}(s,\,t).
\end{equation}
On the other hand, partial differentiation with respect to $t$ of both sides gives
\begin{align}\nonumber
\dfrac{\sum_{n\geq 1}\dfrac{R(n)^{t}}{n^{s}}\ln R(n)}{\sum_{n\geq 1}\dfrac{R(n)^{t}}{n^{s}}}&=
\sum_{p}\dfrac{p^t\ln{p}(p^s-1)}{(p^s+p^t-1)(p^s-1)}\\ \label{3}
&=\sum_{p}\dfrac{p^{t}}{p^s-1+p^t}\ln{p}\,\,\,:=\,\mathcal{T}(s,\,t).
\end{align}

\vspace{1cm}
\hspace{-0.7cm}\textbf{4 Legitimacy of differentiations}\vspace{0.5cm}\\
Above differentiations are legitimate because the derived series $\mathcal{S}(s,\,t)$ and $\mathcal{T}(s,\,t)$ are convergent is RC. We first show this for $\mathcal{T}(s,\,t)$ by using the inequality $\ln{x}<1-x$ for $x>0$. Since the $n$-th prime $p_{n}$ is greater than $n$ [1], we have
\begin{align*}
\dfrac{p^{t}_{n}\ln p_{n}}{p^{s}_{n}-1+p^{t}_{n}}&=\dfrac{\ln p_{n}}{p_{n}^{s-t}-{p^{-t}_{n}}+1}\\
&<\dfrac{\ln p_{n}}{p_{n}^{s-t}}\\
&<\dfrac{p_{n}-1}{p_{n}^{s-t}}\,\,\,\,(\mbox{use of}\,\ln x<1-x,\,x>0)\\
&<\dfrac{1}{p_{n}^{s-t-1}}\\
&<\dfrac{1}{n^{s-t-1}}\,\,\,\,(\mbox{use of}\,p_{n}>n).
\end{align*}
Summing over all primes $p_{n}$ we therefore get
\[\mathcal{T}(s,\,t)=\sum_{p}\dfrac{p^{t}\ln p}{p^s-1+p^t}=\sum_{n\geq 1}\dfrac{p_{n}^{t}\ln p_{n}}{p^{s}_{n}-1+p^{t}_{n}}<\sum_{n\geq 1}\dfrac{1}{n^{s-t-1}}.\]
which clearly is convergent in RC.\\
The convergence of $\mathcal{S}(s,\,t)$ in CR follows from that of $\mathcal{T}(s,\,t)$ by considering that $\dfrac{p^{s}}{p^{s}-1}<2$ for all primes $p$ and $s>0$. We have namely
\[\mathcal{S}(s,\,t)=\sum_{p}\dfrac{p^s}{p^s-1}\dfrac{p^t}{p^s-1+p^t}\ln{p}
<2\,\,\mathcal{T}(s,\,t),\]
and as $\dfrac{p^s}{p^s-1}>1$, we also have $\mathcal{T}(s,\,t)<\mathcal{S}(s,\,t)$. Combining, we obtain
\[\mathcal{T}(s,\,t)<\mathcal{S}(s,\,t)<2\,\,\mathcal{T}(s,\,t)\]
or
\begin{equation}\label{4}
1<\dfrac{\mathcal{S}(s,\,t)}{\mathcal{T}(s,\,t)}<2,\hspace{0.3cm}\{s,\,t\}\in\mbox{CR}.
\end{equation}
Moreover, this demonstrates that
\[
1<\underset{\{s,\,t\}\in\mbox{RC}}{\underline{L}}\,\,\,\dfrac{\mathcal{S}(s,\,t)}{\mathcal{T}(s,\,t)}
<\underset{\{s,\,t\}\in\mbox{RC}}{\overline{L}}\,\,\,\dfrac{\mathcal{S}(s,\,t)}{\mathcal{T}(s,\,t)}<2.
\]
We shall not use this inequality here but it looks it is important for more detailed investigations.The more so, as the technique we used so far applies \textbf{as is} to any positive multiplicative function $M(n)$ (subject to convergence questions). Indeed, in such a case we again have
\[
\mathcal{T}_{M(n)}(s,\,t)<\mathcal{S}_{M(n)}(s,\,t)<2\hspace{0.3cm}\mathcal{T}_{M(n)}(s,\,t),
\]
where
\[
\mathcal{S}_{M(n)}(s,\,t)=\sum_{p}\dfrac{p^s}{p^s-1}\dfrac{M(p)^{t}}{p^s-1+M(p)^{t}}\ln{p}
\]
and
\[
\mathcal{T}_{M(n)}(s,\,t)=\sum_{p}\dfrac{M(p)^{t}}{p^s-1+M(p)^{t}}\ln{M(p)}
\]
are the functions corresponding to the functions $\mathcal{S}(s,\,t)$ and $\mathcal{T}(s,\,t)$ if $M(n)=R(n)$.
\vspace{1cm}\\
\hspace{-0.7cm}\textbf{5 The identity\,\, $\sum_{n\geq{1}}\dfrac{R(n)^{t}}{n^s}\ln\dfrac{R(n)^{\mathcal{S}(s,\,t)}}{n^{\mathcal{T}(s,\,t)}}=0$}\vspace{0.5cm}\\
From \eqref{2} and \eqref{3} of section $3$ and since $\sum_{n\geq 1}\dfrac{R(n)^{t}}{n^s}$ is not zero we get (writing henceforth $\mathcal{S}$ for $\mathcal{S}(s,\,t)$ and $\mathcal{T}$ for $\mathcal{T}(s,\,t)$)
\[\mathcal{S}\sum_{n\geq 1}\dfrac{R(n)^{t}}{n^s}\ln R(n)=\mathcal{T}\sum_{n\geq 1}\dfrac{R(n)^{t}}{n^s}\ln n
\]
or equivalently
\begin{equation}\label{5}
\sum_{n\geq 1}\dfrac{R(n)^{t}}{n^s}\ln\dfrac{R(n)^{\mathcal{S}}}{n^{\mathcal{T}}}=0.
\end{equation}
\vspace{1cm}\\
\hspace{-0.7cm}\textbf{6 Interpretations of\,\, $\sum_{n\geq{1}}\dfrac{R(n)^{t}}{n^s}\ln\dfrac{R(n)^{\mathcal{S}}}{n^{\mathcal{T}}}=0$}\vspace{0.5cm}\\
We now take as point of reference the ratio ${\mathcal{S}}/{\mathcal{T}}$ and split \eqref{5} as follows
\[\sum_{n<R(n)^{\mathcal{S}/\mathcal{T}}}\dfrac{R(n)^{t}}{n^s}\ln\dfrac{R(n)^{\mathcal{S}}}{n^{\mathcal{T}}}
+\sum_{n=R(n)^{{\mathcal{S}}/{\mathcal{T}}}}\dfrac{R(n)^{t}}{n^s}\ln\dfrac{R(n)^{\mathcal{S}}}{n^{\mathcal{T}}}
+\sum_{n>R(n)^{{\mathcal{S}}/{\mathcal{T}}}}\dfrac{R(n)^{t}}{n^s}\ln\dfrac{R(n)^{\mathcal{S}}}{n^{\mathcal{T}}}=0.
\]
The second summand is zero because
\[
\sum_{n=R(n)^{{\mathcal{S}}/{\mathcal{T}}}}\dfrac{R(n)^{t}}{n^s}\ln\dfrac{R(n)^{\mathcal{S}}}{n^{\mathcal{T}}}
=\sum_{n=R(n)^{{\mathcal{S}}/{\mathcal{T}}}}\dfrac{n^{(\mathcal{T}/{\mathcal{S}})t}}{n^s}\ln\dfrac{n^{\mathcal{T}}}{n^{\mathcal{T}}}
=\sum_{n=R(n)^{{\mathcal{S}}/{\mathcal{T}}}}\dfrac{1}{n^{\mathcal{S}-(\mathcal{T}/{\mathcal{S}})t}}\ln{1}=0.
\]
As a result we therefore obtain
\[
\sum_{n<R(n)^{{\mathcal{S}}/{\mathcal{T}}}}\dfrac{R(n)^{t}}{n^s}\ln\dfrac{R(n)^{\mathcal{S}}}{n^{\mathcal{T}}}
+\sum_{n>R(n)^{{\mathcal{S}}/{\mathcal{T}}}}\dfrac{R(n)^{t}}{n^s}\ln\dfrac{R(n)^{\mathcal{S}}}{n^{\mathcal{T}}}=0
\]
or
\begin{equation}\label{6}
\sum_{n<R(n)^{{\mathcal{S}}/{\mathcal{T}}}}\dfrac{R(n)^{t}}{n^s}\ln\dfrac{R(n)^{\mathcal{S}}}{n^{\mathcal{T}}}
=\sum_{n>R(n)^{{\mathcal{S}}/{\mathcal{T}}}}\dfrac{R(n)^{t}}{n^s}\ln\dfrac{n^{\mathcal{T}}}{R(n)^{\mathcal{S}}}.
\end{equation}
This identity in $s$ and $t$, which is a different interpretation of \eqref{5}, is fundamental for the sequel.\vspace{0.4cm}\\

\hspace{-0.7cm}\textbf{6 The connection with the $abc$-conjecture}\vspace{0.5cm}\\
The deeper meaning of the identity \eqref{6} of the previous section is that both series \[\sum_{n<R(n)^{{\mathcal{S}}/{\mathcal{T}}}}\dfrac{R(n)^{t}}{n^s}\ln\dfrac{R(n)^{\mathcal{S}}}{n^{\mathcal{T}}}\]
and \[\sum_{n>R(n)^{{\mathcal{S}}/{\mathcal{T}}}}\dfrac{R(n)^{t}}{n^s}\ln\dfrac{n^{\mathcal{T}}}{R(n)^{\mathcal{S}}}\] are not \textbf{empty}, as otherwise this would contradict \eqref{4} of section $4$, unless $n^{\mathcal{T}}=R(n)^{\mathcal{S}}$ identically for \textbf{all} $s$ and $t$ within CR. But this, however, is impossible as shown by the case of squarefree numbers for which $R(n)=n$ would give $n^{\mathcal{T}(s,\,t)}=n^{\mathcal{S}(s,\,t)}$, clearly an absurdity as $\mathcal{T}(s,\,t)<\mathcal{S}(s,\,t)<2\,\,\mathcal{T}(s,\,t)$. From these facts we deduce\vspace{0.3cm}\\
\textbf{Theorem 2}.\hspace{0.5cm}
For all coprime integers $a,\,b$ and $c=a+b$ satisfying $c<R(c)^{{\mathcal{S}}/{\mathcal{T}}}$ we have
\[
a+b<R(abc)^{2}
\]
i.e. Bombieri's $abc$-conjecture/question [3].\vspace{1cm}\\
\hspace{-0.5cm}\textbf{Proof}. Since $c$ can be written in $\dfrac{\phi(c)}{2}$ different ways as a sum of two coprime integers $a$ and $b$ [4], [5] and as by assumption $c<R(c)^{{\mathcal{S}}/{\mathcal{T}}}$ it results that $c$ occurs in the sum
\[
\sum_{n<R(n)^{{\mathcal{S}}/{\mathcal{T}}}}\dfrac{R(n)^{t}}{n^s}\ln\dfrac{R(n)^{\mathcal{S}}}{n^{\mathcal{T}}}.
\]
Consequently, we have
\[
a+b=c<R(c)^{{\mathcal{S}}/{\mathcal{T}}}<R(c)^{{\mathcal{S}}/{\mathcal{T}}}R(ab)^{{\mathcal{S}}/{\mathcal{T}}}<R(abc)^{{\mathcal{S}}/{\mathcal{T}}},
\]
which, because of \eqref{4} of section 4, gives
\[
a+b<R(abc)^{2}.
\]
Q.E.D.
\vspace{1cm}\\
\hspace{-0.7cm}\textbf{Note}. Elevating a positive multiplicative function $M(n)$ to the power $t>0$ in $\sum_{n\geq 1}\dfrac{M(n)^{t}}{n^s}$ and forming the functions $\mathcal{S}_{M}(s,\,t)$ and $\mathcal{T}_{M}(s,\,t)$ is an effective technique but raises also problems and questions depending on the cases examined. For example:\vspace{0.5cm}\\
Regarding our worked out example $M(n)=R(n)$
\vspace{0.5cm}\\
$1.$ Is there a deeper meaning that the exponent $2$ in Bombieri's $abc$-conjecture/question coincides with the upper bound of $\dfrac{\mathcal{S}(s,\,t)}{\mathcal{T}(s,\,t)}$ which is also $2$ ?\vspace{0.5cm}\\
$2.$ Is it true that identically in $s$ and $t$
\[
\underset{\{s,\,t\}\in\mbox{RC}}{\underline{L}}\,\,\,\dfrac{\mathcal{S}(s,\,t)}{\mathcal{T}(s,\,t)}\,\,\,\,
\boldsymbol{=}\underset{\{s,\,t\}\in\mbox{RC}}{\overline{L}}\,\,\,\dfrac{\mathcal{S}(s,\,t)}{\mathcal{T}(s,\,t)}<\,\,2,
\]
or is this only true for some specific constant values of $s$ and $t$ ? If so, would this be a proof is the $abc$-conjecture ?\vspace{0.5cm}\\
$3.$ What is the surface $\dfrac{\mathcal{S}(s,\,t)}{\mathcal{T}(s,\,t)}$ like in CR ? Is analytic continuation into the complex plane of $s$ feasible ? Would this give a functional equation as is the case for Riemann's $Zeta$ function ?\vspace{0.5cm}\vspace{0.5cm}\\
Regarding general $M(n)$'s\vspace{0.5cm}\\
$4.$ In what cases are the sums $\mathcal{S}_{M(n)}(s,\,t)$ and $\mathcal{T}_{M(n)}(s,\,t)$ amenable in the sense that we get closed formulas ? \vspace{0.5cm}\\
$5.$ The set of multiplicative functions has a rich structure. Maybe a detailed investigation, in the same way we did for $R(n)$, of the principal $M(n)$'s researched in the Analytic Theory of Numbers, may reveal unsuspected truths !
\newpage
\hspace{-0.7cm}\textbf{References}\vspace{0.5cm}\\
$[1]$ {Hardy G.H. and Wright G.M., An Introduction to the Theory of Numbers, Fourth Edition, Oxford University Press, 1968, pp. 249-253.}\\ \\
$[2]$ {Chandrasekharan K., Introduction to Analytic Number Theory, Springer Verlag Berlin Heildeberg New York, 1968, pp. 76-77.}\\ \\
$[3]$ {Bombieri Enrico, Forty Years of Effective Results in Diophantine Theory, in W$\ddot{u}$stholtz, Gisbert, A Panorama of Number Theory, Cambridge University Press, 2002, p.206.}\\ \\
$[4]$ {Petridi Constantin M., A strong $``$\,$abc$-conjecture$"$ for certain partitions $a+b$ of $c$,\\ arXiv:math/0511224v3 [math.NT] 1 Mar 2006.}\\ \\
$[5]$ {Petridi Constantin M., The number of equations $c=a+b$ satisfying the $abc$-conjecture, arXiv:0904.1935v1 [math.NT] 13 Apr 2009.}\\ \\

\end{document}